\documentclass[11pt]{article}
\usepackage{amsmath}
\usepackage{amssymb}
\usepackage{amsfonts}
\usepackage{eucal}
\newtheorem{theorem}{Theorem}[section]

\newtheorem{corollary}[theorem]{Corollary}

\newtheorem{remark}[theorem]{Remark}
\newtheorem{example}[theorem]{Example}
\begin{document}
\title{A general correlation inequality and the Almost Sure Local Limit Theorem for random sequences in
the domain of attraction of a stable law
\thanks{The financial supports of the Research Grant PRIN 2008
\emph{Probability and Finance}, of the INDAM--GNAMPA and of the Polish National Science Centre, Grant N N201 608740 are
gratefully acknowledged.}}
\author{Rita Giuliano\thanks{Address: Dipartimento di
Matematica "L. Tonelli", Universit\`a di Pisa, Largo Bruno
Pontecorvo 5, I-56127 Pisa, Italy. e-mail:
\texttt{giuliano@dm.unipi.it}} \and Zbigniew S. Szewczak
\thanks{Address: Nicolaus Copernicus University,
Faculty of Mathematics and Computer Science, ul. Chopina 12/18
87-100 Toru\'n, Poland. e-mail: \texttt{zssz@mat.uni.torun.pl}}}
\date{}
\maketitle
\begin{abstract}

\noindent In the present paper we obtain a new correlation
inequality and use it for the purpose of exten\-ding the theory of
the Almost Sure Local Limit Theorem to the case of lattice random
sequences in the domain of attraction of a stable law.
In particular, we prove ASLLT in the case of the normal domain of attraction of $\alpha$--stable law, $\alpha\in(1,2)$. 


\end{abstract}

\section{Introduction}
In the recent paper \cite{W}, the author proves a correlation inequality and an Almost Sure Local Limit Theorem (ASLLT) for i.i.d. square integrable random variables taking values in a lattice.  The sequence of partial sums of such variables are of course in the domain of attraction of the normal law, which is stable of order $\alpha =2$.

\smallskip
\noindent
The aim of the present paper is to give an analogous correlation inequality (Theorem \ref{cov}) for the more general case of random sequences in the domain of attraction of a stable law of order $\alpha \leqslant 2$ and to apply it for the purpose of extending the theory of ASLLT. Notice that in our situation the summands need not be square integrable. Our correlation inequality turns out to be of the typical form needed in the theory of Almost Sure (Central and Local) Limit Theorems (see Corollary \ref{cor} and Remark \ref{primorem}). Our work is based on a careful use of the form of the characteristic function, and is completely different from the one used in \cite{W} (Mc Donald's method of extraction of the Bernoulli part of a random variable).

\bigskip\noindent
{\bf Acknowledgement.} We are grateful to an anonymous referee whose suggestions have led to a substantial improvement of the presentation.
\section{The assumptions and some preliminaries}

In this paper we shall be concerned with a sequence of i.i.d. random
variables $(X_n)_{n \geq 1}$ such that their common distribution $F$ is in the domain
of attraction of $G$, where $G$ is a stable  distribution with exponent $\alpha$ ($0<\alpha \leqslant
2$, $\alpha \not =1$). This means that, for a suitable choice of constants $a_n$ and
$b_n$, the distribution of
$$T_n := \frac {X_1 + \cdots X_n -a_n }{b_n}$$
converges weakly to $G$. It is well known (see \cite {IL}, p. 46)
that in such a case we have $b_n = L(n) n^{1/\alpha}$, where $L$ is
slowly varying in Karamata's sense. For $\alpha > 1$ we shall assume that $X_1$ is centered; by Remark 2 p. 402 of \cite {AD}, this implies that $a_n =0$, for every $\alpha$.

\bigskip
\noindent
 We shall suppose that $X_1$ takes values in the lattice $\mathcal{L}(a,d)= \{a+ kd, \, k \in \mathbb{Z}\}$ where $d$ is the maximal span of the distribution; hence  $S_n := X_1 + \cdots X_n $ takes values in the lattice $\mathcal{L}(na,d)= \{na+ kd, \, k \in \mathbb{Z}\}.$

\bigskip
\noindent For every $n$, let $\kappa_n$ be a number of the form $na
+ kd$ and let
$$\lim_{n \to \infty}{{\kappa_n}\over{b_n}}= \kappa.$$

\bigskip
\noindent Observe that Theorem 4.2.1 p. 121 in \cite{IL} implies that
\begin{equation}\label{qq}\sup_n\big\{\sup_k b_n P(S_n= k)\big\}= C < \infty.
\end{equation}
\bigskip
\noindent Throughout this paper we assume that
\begin{equation}
\label{e2}
x^\alpha P(X>x) = \big(c_1 +o(1)\big)l(x); \qquad x^\alpha
P(X\leqslant -x) = \big(c_2 +o(1)\big)l(x),\quad
\alpha\in(0,2],\end{equation} where $l$ is slowly varying as $x
\to \infty$ and $c_1$ and $c_2$ are two suitable non--negative
constants, $c_1+c_2>0$, related to the stable distribution $G$.

\bigskip\noindent
Let $\phi$  be the characteristic function of $F$.
By \cite{AD}, Theorem 1, for $\alpha \not =1$ it  has the form
\begin{equation}\label{funzionecaratteristica}
\phi(t) = \exp\left\{-c|t|^\alpha h\Big(\frac{1}{|t|}\Big)(1-i\beta {\rm sign}(t)  \tan\frac{\pi\alpha}{2})+ o\left(|t|^\alpha h\Big(\frac{1}{|t|}\Big)\right)\right\},
\end{equation}
where $c= \Gamma(1-\alpha)(c_1+c_2) \cos \frac{\pi\alpha}{2}>0$ and $\beta =\frac{c_1-c_2}{c_1+c_2}\in[-1,1]$
are two constants and $h(x)=l(x)$ if $\alpha\in(0,2)$ and $c=\frac{1}{2},$ $\beta=0,$ $h(x)=E[X^2 1_{\{|X|\leqslant x\}}]$
if $\alpha=2$.
This formula implies that
\begin{align}\label{logaritmo}
\log \big|\phi (t)\big|& = \mathfrak{Re}\big(\log \phi (t)\big)= -c|t|^\alpha h\Big(\frac{1}{|t|}\Big)\big(1+o(1)\big),
\\ \arg\big( \phi (t)\big)& =\mathfrak{Im}\big(\log \phi (t)\big)= -c|t|^\alpha h\Big(\frac{1}{|t|}\Big)\big(-\beta{\rm sign}(t)\tan \frac {\pi \alpha}{2}  +o(1)\big)\nonumber
 \end{align}
hence
\begin{align} \label{argomentosulogaritmo}\lim_{t \to 0}\Big|\frac {\arg \big(\phi(t)\big)}{\log|\phi(t)|}\Big|=\Big|\beta \tan \frac {\pi \alpha}{2}\Big|.
 \end{align}

 \noindent
 We notice that $L(n)= h^{\frac{1}{\alpha}}(b_n)$ for $\alpha\in(0,2)$ (by Remark 2 p. 402 in \cite{AD}), while
  $L(n)=\sqrt{E\big[X^21_{\{|X|\leqslant b_n\}}\big]}$ for $\alpha=2$.

 \begin{remark}\label{controesempio}\rm
 For the case $\alpha=2$ we need that $x \mapsto x^2 P(|X|>x)$
 is a slowly varying function, a stronger assumption than the slow variation of $x \mapsto E[X^2 1_{\{|X|\leqslant x\}}]$ (which in turn is equivalent to the CLT, see Corollary 1 p. 578 in \cite{F}). To see this, consider the following distribution:

 $$P(X=n) = \frac{C}{n^2 2^n}, \quad n\geq 1
 , \qquad C= \sum_{k \geqslant 1}\frac{1}{k^2 2^k}.$$
 It is easy to check that in this case $x \mapsto x^2 P(|X|>x)$
 is not slowly varying.

 \bigskip
 \noindent

 \end{remark}
\begin{remark}\label{equivalenza} \rm Let $\widetilde h\sim h$ as $x \to + \infty$. Then, by \eqref{logaritmo},
\begin{equation*}
 \log \big|\phi (t)\big| = -c|t|^\alpha h\Big(\frac{1}{|t|}\Big)\big(1+o(1)\big)=  -c|t|^\alpha \widetilde h\Big(\frac{1}{|t|}\Big)\cdot \frac{h\Big(\frac{1}{|t|}\Big)}{\widetilde h\Big(\frac{1}{|t|}\Big)}\big(1+o(1)\big) = -c|t|^\alpha \widetilde h\Big(\frac{1}{|t|}\Big)\big(1+o(1)\big).
\end{equation*}
This means that $h$ is unique up to equivalence; thus, by Theorem 1.3.3. p. 14 of \cite{BGT} we can assume that $h$ is continuous (even $C^\infty$) on $[a, \infty)$ for some $a>0$.

\smallskip
\noindent
An analogous observation is in force for $\arg\big(\phi(t)\big)$.

\end{remark}
\begin{remark}
 \rm Thus we deal with a subclass of strictly stable distributions.
Denoting by $\psi$ the characteristic function of $G$, we know from \cite{Z},
Theorem C.4 on p.17 that  $\log \psi$ (for strictly
stable distributions) has the form
$$\log \psi(t)= - c |t|^\alpha \exp \{-i\left(\frac{\pi}{2}\right) \theta
\alpha\,{\rm sign} (t)\},$$
where $|\theta| \leqslant \min \{1, \frac{2}{\alpha}-1\}$ and $c >0$. For $\alpha=1$
and $|\theta|=1$ we get degenerate
distribution and in this case we say that $X_n$ is relatively stable (see e.g.
\cite{SZ}). Almost sure variant of
relative stability for dependent strictly stationary sequences will be discussed
elsewhere.\end{remark}

\bigskip
\noindent
Let $\delta > -1$ and $p>0$ two given numbers; we shall use the equality
\begin{equation}\label{formula}
\int_0^{+ \infty} t^\delta e^{-p t^\alpha} dt = \frac{\Gamma(\frac{\delta +1}{\alpha})}{\alpha}\cdot \frac{1}{p^{\frac{\delta +1}{\alpha}}}= C  \cdot\frac{1}{p^{\frac{\delta +1}{\alpha}}}.
\end{equation}
\smallskip
\noindent
In what follows, with the symbols $C$, $c$ and so on we shall mean
positive constants the value of which may change from case to case.

\bigskip

\bigskip
\section{The correlation inequality}

We assume that $(X_n)_{n \geqslant 1}$ is a sequence of i.i.d. random variables verifying the following conditions:
(\ref{e2}), $\alpha \not =1$ and
$\mu=E[X_1]=0$ when $\alpha >1$. Recall that the norming constant are $a_n=0$ and $b_n= L(n) n^{1/\alpha}$ with $L$ slowly varying. With no loss of generality, we shall assume throughout that
$d=1$.
\bigskip
\noindent
\begin{theorem}\label{cov} (i) In the above setting we have\begin {align*} &b_mb_n
\Big|P(S_m=\kappa_m,S_n=\kappa_n )- P(S_m=\kappa_m
)P(S_n=\kappa_n)\Big|\\&\leqslant C\Big\{\Big(\frac
{n}{n-m}\Big)^{1/\alpha}\frac {L(n)}{L(n-m)}+1\Big\}.\end{align*}
(ii) \underline {In addition} to the previous hypotheses assume that the function $h$ appearing in \eqref{e2}  and \eqref{funzionecaratteristica} verifies
$$\liminf_{x\to \infty}h(x)=:\ell>0.$$

\bigskip
\noindent
Then there exists $\epsilon >0$ such that, putting
$$M(x) = \sup_{\frac{1}{\epsilon}\leqslant y \leqslant x}h(y) , \qquad x \geqslant  \frac{1}{\epsilon},$$we have $M(x) < \infty$ for every $x$ and
\begin {align}\label{secondacorrelazione}\nonumber & b_mb_n
\Big|P(S_m=\kappa_m,S_n=\kappa_n )- P(S_m=\kappa_m
)P(S_n=\kappa_n)\Big|\\&\leqslant CL(n)\left\{n^{1/\alpha}\Big(\frac {1}{e^{(n-m)c} } + \frac {1}{e^{nc}
}\Big)+ \frac{\frac{m}{n}}{\Big(1- \frac{m}{n}\Big)^{1+ \frac{1}{\alpha}}}\Big(1+ M(n^{1+\frac{1}{\alpha}})\Big) + \frac{\big(\frac{m}{n}\big)^\frac{\eta}{\alpha} L^\eta(m)}{(1-\frac{m}{n})^\frac{\eta +1}{\alpha}}  \right\}\end {align}
for every pair $(m,n)$ of integers, with $m \geqslant 1 $, $n > m+ \epsilon ^{-\frac{\alpha}{\alpha +1}}$, and for every $\eta \in (0,1]$.

\end{theorem}

\begin{remark}\label{nuovo}\rm
If $h$ is ultimately increasing, then condition (ii) of Theorem \ref{cov} is automatically satisfied. A quick look at the proof (see below) shows that if $h$ is increasing and continuous, the inequality \eqref{secondacorrelazione} holds for $1 \leqslant m<n$.
\end{remark}
\emph{Proof of Theorem \ref{cov}.}

\bigskip
\noindent (i) We write
\begin{align*}  &b_mb_n
\Big|P(S_m=\kappa_m,S_n=\kappa_n )- P(S_m=\kappa_m
)P(S_n=\kappa_n)\Big|\\&= \Big\{b_m P(S_m=\kappa_m)\Big\} \cdot
\Big\{b_n \Big|P(S_{n-m}=\kappa_{n}-\kappa_m )-
P(S_n=\kappa_n)\Big|\Big\}\\&\leqslant C\cdot b_n
\Big|P(S_{n-m}=\kappa_{n}-\kappa_m )- P(S_n=\kappa_n)\Big|,
\end{align*} by \eqref{qq}. Inequality (i) follows since
\begin{align*}&b_n \Big|P(S_{n-m}=\kappa_{n}-\kappa_m )-
P(S_n=\kappa_n)\Big|\leqslant b_n \Big(P(S_{n-m}=\kappa_{n}-\kappa_m
)+P(S_n=\kappa_n)\Big)\\&= \Big(\frac {b_n}{b_{n-m}}\cdot
b_{n-m}P(S_{n-m}=\kappa_{n}-\kappa_m )+b_nP(S_n=\kappa_n)\Big)\leqslant C
\Big(\frac {b_n}{b_{n-m}}+1\Big),\end{align*} by \eqref{qq} again.

\bigskip
\noindent (ii) Let $\phi$ be the characteristic function of $F$. By
the inversion formula (see Theorem 4, p. 511 of \cite{F}) we can
write
\begin{align*}\label{bb}& b_n
\Big|P(S_{n-m}=\kappa_{n}-\kappa_m )- P(S_n=\kappa_n)\Big|= \frac
{b_n}{2\pi}\Big|\int_{-\pi}^\pi\Big\{e^{-it
(\kappa_n-\kappa_m)}\phi^{n-m}(t)- e^{-it
\kappa_n}\phi^{n}(t)\Big\}\,dt\Big|\\& \leqslant Cb_n
\int_{-\pi}^\pi\Big|e^{it
\kappa_m}\phi^{n-m}(t)-\phi^{n}(t)\Big|\,dt.
\end{align*}

\bigskip
\noindent
Recall the expression \eqref{funzionecaratteristica} of $\phi$ where,  by Remark \ref{equivalenza}, on can choose $h$ continuous on $[a, \infty)$ for some $a>0.$ The additional assumption $\ell >0$ allows to take $\epsilon \in\big(0, \frac{1}{a}\big]$ such that, for $|t|< \epsilon$ we have
\begin{equation}\label{minimolimite}
h\Big(\frac{1}{|t|}\Big) >\frac{A}{2}>0.
\end{equation}
 We write
$$\int_{-\pi}^\pi\Big|e^{it
\kappa_m}\phi^{n-m}(t)-\phi^{n}(t)\Big|\,dt= \int_{|t|<\epsilon} +
\int_{\epsilon< |t|< \pi}  = I_1+I_2.$$ Now
$$I_2\leqslant  \int_{\epsilon< |t|< \pi}|\phi(t)|^{n-m}\,dt + \int_{\epsilon< |t|< \pi}|\phi(t)|^{n}\,dt. $$
Since $d=1$, by Theorem 1.4.2 p. 27 of \cite{IL} we have $|\phi(t)|<1$ for
$0<|t|< 2\pi$. Hence a constant $c >0$ exists such that, for
$\epsilon < |t| < \pi$ we have $|\phi(t)|< e^{-c}$, which gives $$
\int_{\epsilon< |t|< \pi}|\phi(t)|^{n}\,dt < 2\pi e^{-nc}; \qquad
\int_{\epsilon< |t|< \pi}|\phi(t)|^{n-m}\,dt < 2\pi e^{-(n-m)c},$$
so that
\begin{equation}\label {I2}I_2\leqslant  C \Big(e^{-(n-m)c} +
e^{-nc}\Big).\end{equation} Now we evaluate $I_1$,
\begin{align}\label{stimadiIuno}\nonumber
&I_1 =\int_{-\epsilon}^\epsilon\Big|e^{it
\kappa_m}\phi^{n-m}(t)-\phi^{n}(t)\Big|\,dt\\& \nonumber\leqslant \int_{-\epsilon}^\epsilon\Big|e^{it
\kappa_m}\phi^{n-m}(t)-e^{it
\kappa_m}\phi^{n}(t)\Big|\,dt+\int_{-\epsilon}^\epsilon\Big|e^{it
\kappa_m}\phi^{n}(t)-\phi^{n}(t)\Big|\,dt\\&=\int_{-\epsilon}^\epsilon\Big|\phi^{n-m}(t)-\phi^{n}(t)\Big|\,dt+\int_{-\epsilon}^\epsilon\Big|e^{it
\kappa_m}-1\Big|\cdot\big|\phi(t)\big|^{n}\,dt.
\end{align}
\bigskip \noindent For a complex number $A=\varrho {\rm e}^{i \theta}$ we have
\begin {align*}&\big|A^x- A^y\big| %= \big|\rho^xe^{i\theta x}-
=\Big\{\big( \varrho^{x}- \varrho^{y}\big)^2+ 2
\varrho^{x+y}\big(1-\cos\theta(x-y)\big)\Big\}^{1/2}\leqslant \big|\varrho^{x}-
\varrho^{y}\big|+\Big\{ 2
\varrho^{x+y}\big(1-\cos\theta(x-y)\big)\Big\}^{1/2}\\& \leqslant
\big|\varrho^{x}- \varrho^{y}\big|+\varrho^{\frac{x+y}{2}}
|\theta|\big|x-y\big|.
\end {align*} Applying with $A= \phi(t)$, $x= n-m$ and $y=n$ we get
\begin{equation}\label{boiade}\int_{|t|<\epsilon}\Big|\phi^{n-m}(t)-\phi^{n}(t)\Big|\,dt \leqslant
C\Big(
\int_{|t|<\epsilon}\Big||\phi(t)|^{n-m}-|\phi(t)|^{n}\Big|\,dt
+m\int_{|t|<\epsilon}\Big|\arg
\phi(t)\Big|\cdot\big|\phi(t)\big|^{n-(m/2)}\,dt\Big).\end{equation}
Applying Lagrange Theorem to the first summand we find that for a
suitable $\xi \in (n-m,n)$ we have, for every $\delta < \epsilon$
\begin{align}& \label{r}\nonumber
\int_{|t|<\epsilon}\Big||\phi(t)|^{n-m}-|\phi(t)|^{n}\Big|\,dt \leqslant \int_{|t|< \delta}2\, dt+
m\int_{\delta <|t|<\epsilon}\Big|\frac{d}{dx}\{|\phi(t)|^x\}
\Big|\Bigg|_{x=\xi}\,dt\\  \nonumber &=4 \delta +m
\int_{\delta <|t|<\epsilon}\big|\log |\phi(t)|\big|\cdot |\phi(t)|
^\xi\,dt\leqslant 4 \delta +m \int_{\delta <|t|<\epsilon}\big|\log |\phi(t)|\big|\cdot
|\phi(t)| ^{(n-m)}\,dt\\   &\leqslant 4 \delta + C_1 m
\int_{\delta <|t|<\epsilon}|t|^\alpha h\Big(\frac{1}{|t|}\Big)\cdot
e ^{-C_2(n-m)|t|^\alpha h\big(\frac{1}{|t|}\big)}\,dt,
\end{align}
using the relation \eqref{logaritmo}. By reporting the inequality \eqref{minimolimite} into \eqref{r}, and recalling that $h$ is continuous, hence bounded on $[\delta, \epsilon]$, we obtain
\begin{align*}&
\int_{|t|<\epsilon}\Big||\phi(t)|^{n-m}-|\phi(t)|^{n}\Big|\,dt \leqslant 4\delta +C_1 m M\Big(\frac{1}{\delta}\Big)
\int_{\delta<|t|<\epsilon}|t|^{\alpha}  \cdot
e ^{-C_2(n-m)|t|^{\alpha}} \,dt \\&\leqslant C\Big\{\delta + \frac{m}{(n-m)^{1+ \frac{1}{\alpha}}}\cdot M\Big(\frac{1}{\delta}\Big)\Big\},\end{align*}
by \eqref{formula}, for any $\delta < \epsilon$. Taking
$\delta = \frac{1}{(n-m)^{1+ \frac{1}{\alpha}}}$,
we get
\begin{align}\label{basic1}
\nonumber &\int_{|t|<\epsilon}\Big||\phi(t)|^{n-m}-|\phi(t)|^{n}\Big|\,dt \leqslant C\Big\{\frac{1}{(n-m)^{1+ \frac{1}{\alpha}}}+ \frac{m}{(n-m)^{1+ \frac{1}{\alpha}}}\cdot M\big((n-m)^{1+ \frac{1}{\alpha}}\big)\Big\}\\& \leqslant C\frac{m}{(n-m)^{1+ \frac{1}{\alpha}}}\Big(1+ M\big(n^{1+ \frac{1}{\alpha}}\big)\Big),
\end{align}
$M$ being non--decreasing.

\bigskip
\noindent
For the second
summand in \eqref{boiade} we can proceed as follows: by \eqref{argomentosulogaritmo},
$$\big|\arg \phi(t)\big|\leqslant C \big|\log |\phi(t)|\big|, \qquad \forall \, t.$$
Hence, arguing as before
\begin {align}\label{Y}
&\nonumber m\int_{|t|<\epsilon}\big|\arg\phi(t)\big|\cdot\big|\phi(t)\big|^{n-(m/2)}\,dt
 \nonumber\leqslant Cm
\int_{|t|<\epsilon}\big|\log\phi(t)\big|\cdot\big|\phi(t)\big|^{n-(m/2)}\,dt\\&  \leqslant
  C \frac{m}{\big(n-
\frac{m}{2}\big)^{1+\frac{ 1}{\alpha }}}\Big\{1+ M\Big(\big(n-\frac{m}{2}\big)^{1+ \frac{1}{\alpha}}\Big)\Big\}\leqslant C \frac{m}{(n-m)^{1+\frac{1}{\alpha }}}\Big(1+ M\big(n^{1+ \frac{1}{\alpha}}\big)\Big),\end {align}
as before. Thus, by \eqref{boiade} \eqref{basic1}, \eqref{Y} we obtain
\begin{equation}\label{primoaddendo}
\int_{-\epsilon}^\epsilon\Big|\phi^{n-m}(t)-\phi^{n}(t)\Big|\,dt \leqslant C\frac {m}{(n-m)^{1+\frac{1}{\alpha }}}\Big(1+ M\big(n^{1+ \frac{1}{\alpha}}\big)\Big).
\end{equation}
Let's turn to the second summand in \eqref{stimadiIuno}. By the well known inequality$$|e^{it}-1|\leqslant 2^{1-\eta}|t|^\eta , \qquad \forall \eta \in (0,1]$$ (see \cite{L}, p. 200), we have
\begin{align} \label{secondaparte}
&\int_{-\epsilon}^\epsilon\Big|e^{it
\kappa_m}-1\Big|\cdot\big|\phi(t)\big|^{n}\,dt\leqslant |\kappa_m|^\eta 2^{1-\eta} \int_{-\epsilon}^\epsilon\big|t\big|^\eta\cdot \big|\phi(t)\big|^{n-m}\,dt\leqslant C \frac{|\kappa_m|^\eta}{(n-m)^{\frac{\eta +1}{\alpha}}}\leqslant C \frac{m^\frac{\eta}{\alpha}L^\eta(m)}{(n-m)^{\frac{\eta +1}{\alpha}}},
\end{align}
again by \eqref{formula} and the fact that $$\kappa_m \sim \kappa b_m = \kappa L(m) m ^{1/\alpha}, \qquad m \to \infty.$$
Summing the estimates \eqref{I2}, \eqref{primoaddendo} and
\eqref{secondaparte} we get, for every $\eta\in(0,1]$
\begin{align*}
&\int_{-\pi}^{-\pi}\Big|e^{it\kappa_m}\phi^{n-m}(t) - \phi^{n}(t)\Big|\,dt\\& \leqslant C \Bigg\{\Big(e^{-(n-m)c} +
e^{-nc}\Big)+ \frac {m}{(n-m)^{1+\frac{1}{\alpha }}}\Big(1+ M\big(n^{1+ \frac{1}{\alpha}}\big)\Big)+\frac{m^\frac{\eta}{\alpha}L^\eta(m)}{(n-m)^{\frac{\eta +1}{\alpha}}}\Bigg\}.
\end{align*}

\noindent
Multiplying by $b_n =L(n)
n^{1/\alpha}$ gives the conclusion.

\hfill$\Box$

\bigskip
\noindent
\begin{corollary}\label{cor} For large $m$ and $n \geq 2m$, for every $\delta < \frac{1}{\alpha}$ and for every $\eta \in(0,1]$ we have $$b_mb_n
\Big|P(S_m=\kappa_m,S_n=\kappa_n )- P(S_m=\kappa_m
)P(S_n=\kappa_n)\Big| \leqslant C \tilde L(n)\cdot \Big(\frac{m}{n}\Big)^\rho,$$
with $\tilde L(n) =L(n)\big(1+ M(n^{1+\frac{1}{\alpha}})+L^\eta(n)\big)$ and $\rho:= \min\{\eta (\frac{1}{\alpha}- \delta),1\}$.
\end {corollary}

\bigskip
\noindent \emph{Proof of Corollary \ref{cor}.} Let $\epsilon$ be the number identified in (ii) of Theorem \ref{cov} and $c$ the constant appearing in the second member of \eqref{secondacorrelazione}; let $x_0> \epsilon^{-\frac{\alpha}{\alpha +1}}$ be such that $e^{cx}\geq x^{
2/\alpha}$ for $x \geq x_0$. For $m \geq x_0$ we have also $n-m \geqslant
x_0 \geqslant \epsilon^{-\frac{\alpha}{\alpha +1}}$ (since $n-m \geq m$). Then (ii) of Theorem \ref{cov} holds and
\begin{equation}\label{prima}
\frac {n^{1/\alpha}}{e^{nc} }\leqslant \frac {n^{1/\alpha}}{e^{(n-m)c} }\leqslant \frac {n^{1/\alpha}m^{ 1/\alpha}}
{(n-m)^{2/\alpha} }=\frac{\big(\frac{m}{n}\big)^{
1/\alpha}}{\big(1-\frac{m}{n}\big)^{2/\alpha}}\leq 2^{2/\alpha}\cdot
\Big(\frac{m}{n}\Big)^{ 1/\alpha}.
\end{equation}
Moreover

\begin{equation}\label{seconda}
\frac {\frac{m}{n}}{\big(1-\frac {m}{n}\big)^{1+1/\alpha}}
\leqslant
\frac {\frac{m}{n}}{\big(\frac {1}{2}\big)^{1+1/\alpha}}=
(2^{1+\frac{1}{\alpha}})\Big(\frac{m}{n}\Big);\end{equation}
similarly
\begin{equation}\label{secondoaddendo}
\frac{\big(\frac{m}{n}\big)^\frac{\eta}{\alpha} L^\eta(m)}{(1-\frac{m}{n})^\frac{\eta +1}{\alpha}}\leqslant 2^\frac{\eta +1}{\alpha} L^\eta(m)\Big(\frac{m}{n}\Big)^\frac{\eta}{\alpha}. \end{equation}
Recall the well known representation of slowly varying functions (see for instance \cite{BGT}, p.12):
$$L(x) = \gamma(x) \exp \Big\{\int_1 ^x\frac{\varepsilon(t)}{t}\, dt\Big\},$$
where $\gamma(x)\to \gamma$ (a finite constant) and $\varepsilon(x)\to 0$ as $x \to \infty$.

\noindent
We deduce from it that, for every $\delta > 0$, $n \geqslant 2m$ and large $m$ we have
\begin{align*}
&\frac{m^\delta L(m)}{n^\delta L(n)}\leqslant C \exp \Big\{\delta \log m + \int_1 ^m \frac{\varepsilon(t)}{t}\, dt-\delta \log n - \int_1 ^n \frac{\varepsilon(t)}{t}\, dt \Big\}= C \exp\Big \{\delta \log \frac{m}{n} + \int_n ^m \frac{\varepsilon(t)}{t}\, dt \Big\}\\ &=C \exp\Big \{\delta \log \frac{m}{n} - \int_m ^n \frac{\varepsilon(t)}{t}\, dt \Big\} \leqslant C \exp\Big \{\delta \log \frac{m}{n} + \int_m ^n \frac{\delta}{2t}\, dt \Big\}=C \exp\Big \{\delta \log \frac{m}{n}-\frac{\delta}{2} \log \frac{m}{n}\Big\}\\& = C \Big(\frac{m}{n}\Big)^{\frac{\delta}{2} }\leqslant C \Big(\frac{1}{2}\Big)^{\frac{\delta}{2} }.
\end{align*}
It follows that

\begin{equation}\label{maggiorazione}
 L^\eta(m)\Big(\frac{m}{n}\Big)^\frac{\eta}{\alpha}\leqslant C
L^\eta(n)\Big(\frac{m}{n}\Big)^{\eta (\frac{1}{\alpha}- \delta)}.\end{equation}
From \eqref{secondoaddendo} and \eqref{maggiorazione} we obtain
\begin{equation}\label{finale}
\frac{\big(\frac{m}{n}\big)^\frac{\eta}{\alpha} L^\eta(m)}{(1-\frac{m}{n})^\frac{\eta +1}{\alpha}}\leqslant CL^\eta(n)\Big(\frac{m}{n}\Big)^{\eta (\frac{1}{\alpha}- \delta)}.
\end{equation}

\medskip
\noindent Now the desired conclusion follows from \eqref{prima}, \eqref{seconda} and \eqref{finale} and the inequality in (ii) of
Proposition \ref{cov}.

\hfill$\Box$

\begin{remark}\label{primorem}
 \rm Let $(Y_n)_{n \geqslant 1})$ be i.i.d. centered random variables with second moments. It is well known that correlation inequalities of the form
 \begin{equation}\label{formagiusta}
 \Big|Cov(Y_m, Y_n)\Big|\leqslant C \Big(\frac{m}{n}\Big)^\rho
 \end{equation}
for some positive constant $\rho$
are useful tools in order to prove Almost Sure Theorems with logarithmic weights, i.e. statements of the form
\begin{equation}\label{ASLT}
\lim_{N \to \infty}\frac{1}{\log
N}\sum_{n=1}^N \frac{Y_n}{n}=0.
\end{equation}
See for instance \cite{GW}, Theorem (2.9) as a reference.

\noindent
The correlation inequality of Corollary \ref{cor} is similar to \eqref{formagiusta}, but notice that the coefficient $\tilde L$ need not be bounded. See also Remark \ref{secondorem}.
\end{remark}
\begin{remark}\label{remW}
 \rm In the case $\alpha =2$ the correlation inequality of Corollary \ref{cor} furnishes $\rho= \eta\big(\frac{1}{2}-\epsilon\big)<  \frac{1}{2}$, while in \cite{W} the better exponent $\rho =\frac{1}{2}$ is found. Nevertheless, we are able to prove an Almost Sure Local Theorem even with the weaker exponent, as we shall see in the next section.
\end{remark}

\section{Application to the Almost Sure Local Limit Theorem}
In this section we apply our main result to prove a suitable form of the Almost Sure Local Limit Theorem.
Denote by $g$ the $\alpha$--stable density function related to the distribution function $G$.
We point out that here we consider only the case $\alpha > 1$. Precisely
\bigskip
\noindent
\begin{theorem}\label{princ} Let $(X_n)_{n \geqslant 1}$ be
a centered, independent and identically lattice distributed (i.i.l.d.) random sequence with span $d=1$; assume moreover that (\ref{e2}) holds with $\alpha\in (1,2]$
and that there exists $\gamma \in (0,2)$
such that$$\sum_{k=a}^b \frac { L(k)\big\{1+ M\big(k^{1+\frac{1}{\alpha}}\big)+L^\eta (k)\big\}}{k} \leqslant C(\log^\gamma b-
\log^\gamma a),$$
for some $\eta \in (0,1]$. If the  condition (ii) of Theorem \ref{cov} is satisfied, then $$\lim_{N \to \infty}\frac{1}{\log
N}\sum_{n=1}^N \frac{b_n}{n}1_{\{S_n = \kappa_n\}}=
g(\kappa).$$\end{theorem}

\bigskip
\noindent
\begin{example}\rm  Let $h(x) = \log^\sigma x$, with $0<\sigma< \frac{\alpha}{1+\alpha}$. Notice that
\begin{equation}\label{sigma}
\frac{\sigma}{\alpha}< \sigma\wedge(1-\sigma)
\end{equation}
Remark 2 p. 402 in \cite{AD} assures that $$b_n^\alpha= n\log^\sigma b_n.$$
Putting $f(x) = \frac{x^\alpha}{\log ^\sigma x }$ and observing that $f$ is strictly increasing for $x> e^{\frac{\sigma}{\alpha}}$, this means that
\begin{equation}\label{rappresentazione}
L(n)= \frac{b_n}{n^{\frac{1}{\alpha}}}=\frac{f^{-1}(n)}{n^{\frac{1}{\alpha}}},
\end{equation}
for sufficiently large $n$. It is not difficult to check that for sufficiently large $n$ $$L(n)\leqslant \log^\delta n, \qquad \forall \delta >\frac{\sigma}{\alpha}.$$ In fact by
\eqref{rappresentazione} this is equivalent to
\begin{equation*}n \leqslant f\big( n^{\frac{1}{\alpha}}\cdot  \log^\delta n\big)= \frac{n\log^{\alpha\delta} n}{\big(\frac{1}{\alpha }\log n + \delta \log \log n\big)^\sigma},
\end{equation*}
which clearly holds for $\alpha \delta > \sigma$. Thus
\begin{equation*}
 L(n)\big\{1+ M\big(n^{1+\frac{1}{\alpha}}\big)+L^\eta (n)\big\}\leqslant \log^\delta n\Big\{1+ \Big(1+\frac{1}{\alpha}\Big)^\sigma\log^\sigma n+ \log^{\delta\eta} n\Big\}\leqslant C\big(\log n\big)^{2\delta\vee (\delta +\sigma)  }.
\end{equation*}
If $\delta < \sigma$ we have $2\delta\vee (\delta +\sigma)=\delta +\sigma$,
hence
$$\sum_{k=a}^b \frac { L(k)\big\{1+ M\big(k^{1+\frac{1}{\alpha}}\big)+L^\eta (k)\big\}}{k} \leqslant C\sum_{k=a}^b \frac { \big(\log k\big)^{\delta +\sigma }}{k}\leqslant C(\log^\gamma b-
\log^\gamma a)$$
with $\gamma = \delta + \sigma +1 < 2$ for any $\delta \in\Big(\frac{\sigma}{\alpha}, \sigma\wedge(1-\sigma)\Big)$ (see \eqref{sigma}). \end{example}
\begin{remark}\rm \label{secondorem}
\rm Let $Y_n  = b_n \Big(1_{\{S_n=\kappa_n\}}-P(S_n = \kappa_n)\Big)$.
Theorem \ref{princ} states exactly relation \eqref{ASLT}: just observe that
$$\lim_{N \to \infty}\frac{1}{\log
N}\sum_{n=1}^N \frac{b_nP(S_n = \kappa_n)}{n}=g(\kappa),$$
by Theorem 4.2.1 p. 121  of \cite{IL}. Of course, our statement requires an auxiliary hypothesis, due to the fact that in the second member of our correlation inequality we have a supplementary factor, $\tilde L(n)$, which need not be bounded, as observed before (Remark \ref{primorem}).
\end{remark}

\begin{remark}\rm If $L \equiv$ a constant  (i.e. $F$ belongs to the domain of {\it normal}  attraction of $G$, according to the definition on p. 92 of \cite{IL}), then the assumption in Theorem \ref{princ} and condition (ii) of Theorem \ref{cov} are automatically satisfied; thus we get the following nice result
\begin{corollary}\label{carino}
If $(X_n)_{n \geqslant 1}$ is a centered i.i.l.d. random sequence with span $d=1$, (\ref{e2}) holds
with $\alpha\in (1,2)$ and $L \equiv c$, then $$\lim_{N \to \infty}\frac{1}{\log
N}\sum_{n=1}^N \frac{c}{n^{1- \frac{1}{\alpha}}}1_{\{S_n = \kappa_n\}}=
g(\kappa).$$
\end{corollary}
\end{remark}

\bigskip
\noindent \emph{Proof of Theorem \ref{princ}.} We shall denote $\tilde L(n) :=L(n)\big\{1+ M\big(n^{1+\frac{1}{\alpha}}\big)+L^\eta (n)\big\}$ as in Corollary \ref{cor} of the previous section. By Ex. 1.11.4  p. 58 of \cite{BGT}, $M$ is slowly varying, hence the same happens for $n \mapsto M\big(n^{1+\frac{1}{\alpha}}\big)$ and for $\tilde{L}$.  Put
 $$
 Z_n := \sum_{k= 2^{n-1}}^{2^n -1}\frac {Y_k}{k},$$
 where $Y_n$ is defined in Remark \eqref{secondorem}. Following the same argument as in \cite{GW1}, we must prove that
$$\lim_{n \to \infty}\frac {\sum_{i=1}^n Z_i}{n}=0.$$
We shall use the Gaal--Koksma Strong Law of Large Numbers, i.e. (see
\cite{PS}, p. 134); here is the precise statement: \begin{theorem}
\label{GK}Let $(Z_n)_{n\geqslant 1}$ be a sequence of centered random variables
with finite variance. Suppose that there exists a constant $\beta >
0$ such that, for all integers $m \geqslant 0$, $n > 0$, \begin
{equation}\label{bound}
E\Big[\big(\sum_{i=m+1}^{m+n}Z_i\big)^2\Big]\leqslant C\big((m+n)^\beta-
m^\beta\big),\end {equation} for a suitable constant $C$ independent
of $m$ and $n$. Then, for each $\delta > 0$, $$\sum_{ i=1}^n Z_i =
O( n^{\beta/2}(\log n)^{2+\delta} ), \quad P- a.s.$$
\end{theorem}
\begin {remark} \label {ultimate} \rm It is easy to see that Theorem \ref{GK}
 is in force even if the bound
\eqref{bound} holds only for all integers $m \geq h_0$, $n> 0$,
where $h_0$ is an integer strictly greater than 0: just take $Z_i=0$
for $i= 1, 2, \dots, h_0$ and use Theorem \ref{GK}. \end {remark}
We go back to the proof of Theorem \ref{princ}, where we shall
repeatedly use Remark \ref{ultimate} without mentioning it. Since
\begin{equation}\label{init}{\rm \bf E}\Big[\big(\sum_{i=m+1}^{m+n}Z_i\big)^2\Big]=
\sum_{i=m+1}^{m+n}{\rm \bf E }[Z_i^2]+2\sum_{m+1\leqslant i<j\leqslant
m+n}{\rm \bf E }[Z_iZ_j],\end{equation} we bound separately these
two summands. We have first \begin{align} \label{0}&{\rm
\bf E }[Z_i^2]=\sum_{h,k=2^{i-1}}^{2^i-1}\frac {1}{hk}{\rm \bf E
}[Y_hY_k]=\sum_{h=2^{i-1}}^{2^i-1}\frac {1}{h^2}{\rm \bf E
}[Y_h^2]+2\sum_{2^{i-1}\leqslant h<k\leqslant 2^i-1}\frac {1}{hk}{\rm \bf E
}[Y_hY_k].
\end{align}
Now, by \eqref{qq}
\begin {equation*}{\rm \bf E
}[Y_h^2] = b_h^2\Big\{P(S_h = \kappa_h)-P^2(S_h =
\kappa_h)\Big\}\leqslant b_h^2P(S_h = \kappa_h)\leqslant C b_h = C \cdot L(h)
h^{1/\alpha}.\end{equation*} Fix any $\epsilon \in(0, 1-
\frac{1}{\alpha})$ and let $h_0$ be such that $L(t) < t^\epsilon$
for $t \geqslant h_0$. Let $m$ be such that $2^m\geqslant h_0$; for $h \geqslant
2^{i-1}\geq 2^m$ we have from the above that ${\rm \bf E }[Y_h^2]
\leqslant C h^{\epsilon+ (1/\alpha )}$, which gives
\begin {equation} \label{1}\sum_{h=2^{i-1}}^{2^i-1}\frac
{1}{h^2}{\rm \bf E }[Y_h^2] \leqslant C \sum_{h=2^{i-1}}^{2^i-1}\frac
{1}{h^{2-\epsilon -1/\alpha}}\leqslant C\cdot\frac {2^{i}
-2^{i-1}}{(2^{i-1})^{2-\epsilon -1/\alpha}}= \frac
{C}{(2^{i-1})^{1-\epsilon -1/\alpha}}\leqslant C.
\end{equation}
Moreover, by (i) of Theorem \ref{cov}, we have
\begin{align}\label{2}&\sum_{2^{i-1}\leqslant h<k\leqslant 2^i-1}\frac {1}{hk}{\rm \bf E
}[Y_hY_k]\leqslant C \sum_{2^{i-1}\leqslant h<k\leqslant
2^i-1}\frac {1}{hk} \Big(\frac {k}{k-h}\Big)^{1/\alpha}\frac
{L(k)}{L(k-h)}+ C\sum_{2^{i-1}\leqslant h<k\leqslant 2^i-1}\frac {1}{hk}.
\end{align}
Now
 \begin{equation}\label{3}\sum_{2^{i-1}\leqslant h<k\leqslant 2^i-1}\frac {1}{hk}=
\sum_{k=2^{i-1} }^{2^i-1}\frac {1}{k} \sum_{ h=2^{i-1} }^{k-1}\frac
{1}{h}\leqslant \Big(\sum_{k=2^{i-1} }^{2^i-1}\frac {1}{k}\Big)^2\leqslant
C;\end{equation} and
\begin{align}\label{4}& \sum_{2^{i-1}\leqslant h<k\leqslant 2^i-1}\frac {1}{hk}
\Big(\frac {k}{k-h}\Big)^{1/\alpha}\frac {L(k)}{L(k-h)} \leqslant  2\sum_{k=2^{i-1}
}^{2^i-1}\frac {L(k)}{k^{2-1/\alpha}} \sum_{ j=1}^{k-2^{i-1}}\frac
{1}{j^{1/\alpha}}\frac {1}{L(j)} \end{align} (we have used the fact that $\frac {k}{2}< 2^{i-1}\leqslant h$).

\bigskip
\noindent
 Now we
are concerned with the inner sum in the last member of \eqref{4}. The
function
$$t \mapsto U(t)=\frac
{1}{[t]^{1/\alpha}}\frac {1}{L([t])}, \qquad t \geq 1$$ is regularly
varying with exponent $-\frac {1}{\alpha}$; hence, from Theorem 1 p.
281 part (b) of \cite{F} we deduce that, for every $p\geqslant
\frac{1}{\alpha}-1$ $$\frac{k^{p+1}U(k)}{\int_1^kx^pU(x) \, dx}\to
p- \frac{1}{\alpha}+1, \qquad k \to \infty .$$ Since
$$\int_1^{k}x^pU(x) \, dx= \sum_{j=2}^k \int_{j-1}^j x^p U(x) \, dx \geqslant
\sum_{j=2}^{k} (j-1)^p\int_{j-1}^j  U(x) \, dx =\sum_{j=1}^{k-1} j^p U(j), $$ we get
$$\liminf_{k \to \infty}\frac{k^{p+1}U(k)}{\sum_{j=1}^{k-1} j^pU(j)}\geqslant \lim_{k \to \infty} \frac{k^{p+1}U(k)}{\int_1^kx^pU(x)
\, dx}= \Big(p- \frac{1}{\alpha}+1\Big), \qquad k \to \infty .$$
 In particular, for $p=0$ we obtain (\underline{remember that $\frac{1}{\alpha}< 1$})
$$\frac{k}{k^{1/\alpha}}\frac
{1}{L(k)}= k U(k)\geqslant C \sum_{j=1}^{k-1} U(j)= C\sum_{ j=1}^{k-1}\frac
{1}{j^{1/\alpha}}\frac {1}{L(j)},$$ whence
$$\sum_{ j=1}^{k-2^{i-1}}\frac {1}{j^{1/\alpha}}\frac
{1}{L(j)} \leqslant \sum_{ j=1}^{k-1}\frac {1}{j^{1/\alpha}}\frac
{1}{L(j)}\leqslant C \frac{k}{k^{1/\alpha}}\frac {1}{L(k)},$$ and
continuing \eqref{4} we obtain
\begin{equation}\label{5}\sum_{k=2^{i-1} }^{2^i-1}\frac
{L(k)}{k^{2-1/\alpha}} \sum_{ j=1}^{k-2^{i-1}}\frac
{1}{j^{1/\alpha}}\frac {1}{L(j)} \leqslant C\sum_{k=2^{i-1}
}^{2^i-1}\frac {L(k)}{k^{2-1/\alpha}}\frac{k}{k^{1/\alpha}}\frac
{1}{L(k)}=C.\end{equation}
Summarizing , from \eqref{2}, \eqref{3}, \eqref{4} and \eqref{5} we
have found \begin{equation}\label{6}\sum_{2^{i-1}\leqslant h<k\leqslant
2^i-1}\frac {1}{hk}{\rm \bf E }[Y_hY_k]\leqslant C,\end{equation} so that
by \eqref{0}, \eqref{1} and \eqref{6} we get
\begin{equation}\label{sec}{\rm \bf E }[Z_i^2]\leqslant C ;\end{equation}
this implies
\begin{equation}\label{third}\sum_{i=m+1}^{m+n}{\rm \bf E }[Z_i^2]\leq Cn \end{equation}
which bounds the first sum in \eqref{init}. Now we consider the
second one, i.e.
$$\sum_{m+1\leqslant i<j\leqslant m+n}{\rm \bf E }[Z_iZ_j].$$
We start with a bound for the summand ${\rm \bf E }[Z_iZ_j]$ when $j
\geqslant i+2$. In this case we have
$$h \leqslant 2^i \leq 2^{j-2} \leqslant \frac{k}{2}.$$
Let $m$ be such that $2^m > x_0$, where $x_0$ is as in Corollary
\ref{cor}. For $i \geqslant m+1$, the same Corollary assures that
\begin{equation}\label{7}{\rm \bf E }[Z_iZ_j]=
\sum_{h= 2^{i-1}}^{2^i -1}\sum_{k= 2^{j-1}}^{2^j -1}
 \frac {1}{hk}{\rm \bf E }[Y_hY_k]\leqslant C \sum_{h= 2^{i-1}}^{2^i -1}
 \frac {1}{h^{1- \rho}}
 \sum_{k= 2^{j-1}}^{2^j -1}
 \frac {\tilde{L}(k)}{k^{1+ \rho}}.\end{equation}
 The function
 $$V(t) = \frac{\tilde{L}([t])}{[t]^{1+ \rho}}, \qquad t \geq 1$$
 is regularly varying with exponent $-(1 +\rho)$.
  Hence,  by Theorem 1 p. 281 part (a)
of \cite{F}, we have
$$\frac {k^{p+1}V(k)}{\int_{k}^{\infty} x^p V(x)\,dx}\to -p +\rho$$
 if $-p+\rho \geqslant 0$ and
$\int_{k}^{\infty} x^p V(x)\,dx$  is finite. In particular we can
take $p=0$, since $$\int_{1}^{\infty} V(x)\,dx=
\sum_{j=1}^{\infty} \frac{\tilde{L}(j)}{j^{1 + \rho}}< + \infty,$$ and
we obtain
\begin{equation}\label{lim}\frac{\frac{\tilde{L}(k)}{k^{ \rho}}}{\int_{k}^{\infty} V(x)\,dx}
=\frac {kV(k)}{\int_{k}^{\infty} V(x)\,dx}\to \rho.\end{equation} Now
\begin{equation}\label{comp1}\int_k^{\infty}
V(x)\,dx=\int_k^{\infty} \frac{\tilde{L}([x])}{[x]^{ 1+\rho}}
\,dx=\sum_{j=k}^{\infty}\int_j^{j+1} \frac{\tilde{L}([x])}{[x]^{
1+\rho}} \,dx\geqslant\sum_{j=k}^{\infty} \frac{\tilde{L}(j)}{(j+1)^{
1+\rho}}\geqslant
\Big(\frac{1}{2}\Big)^{1+\rho}\sum_{j=k}^{\infty}
\frac{\tilde{L}(j)}{j^{ 1+\rho}}\end{equation} and similarly
 \begin{equation}\label{comp2}\int_k^{\infty} V(x)\,dx\leqslant\sum_{j=k}^{\infty}
\frac{\tilde{L}(j)}{j^{ 1+\rho}}.\end{equation}
 From \eqref{lim}, \eqref{comp1}  and \eqref{comp2} we deduce that there exist two constants $0<C_1< C_2$
such that, for every sufficiently large $k$,
$$C_1 \frac {\tilde{L}(k)}{k^{\rho}} <\sum_{j=k}^{\infty}
 \frac {\tilde{L}(k)}{k^{ 1+\rho}}< C_2 \frac {\tilde{L}(k)}{k^{ \rho}} .$$
Going back to \eqref{7}, we find for sufficiently large $i$
 $${\rm \bf E
}[Z_iZ_j]\leq  2^{i\rho}\Big(C_2
 \frac{\tilde{L}(2^{j-1})}{2^{(j-1)\rho}}- C_1 \frac {\tilde{L}(2^j)}{2^{j\rho}}\Big),$$
 and now, by \eqref{sec}
 \begin{align}& \nonumber \label {fourth}\sum_{m+1\leqslant i<j\leqslant m+n}{\rm \bf E
}[Z_iZ_j]= \sum_{m+3 \leqslant i+2\leqslant j\leqslant m+n}{\rm \bf E }[Z_iZ_j] +
\sum_{ i=m+1}^{m+n-1}{\rm \bf E }[Z_iZ_{i+1}]=
\\ \nonumber & \leqslant
\sum_{j=m+2}^{m+n}\Big( C_2
 \frac{\tilde{L}(2^{j-1})}{2^{(j-1)\rho}}- C_1 \frac {\tilde{L}(2^j)}{2^{j\rho}}
 \Big)\sum_{i=m+1}^{j-1}2^{i\rho}+
\sum_{ i=m+1}^{m+n-1}{\rm \bf E }[Z_i^2]^{1/2}{\rm \bf E
}[Z_{i+1}^2]^{1/2}
\\& \leqslant
C\Big(\sum_{j=m+1}^{m+n-1} \tilde{L}(2^j)+n\Big)\leq C\big\{(m+n)^{\gamma}-
m^{\gamma}+n\}.\end{align}

\bigskip
\noindent Now we insert \eqref {third} and \eqref{fourth} into
\eqref{init} and obtain
$${\rm \bf E}\Big[\big(\sum_{i=m+1}^{m+n}Z_i\big)^2\Big]\leqslant C\big\{(m+n)^{\gamma}-
m^{\gamma}+n\}\leq C\big\{(m+n)^{\gamma\vee1}- m^{\gamma\vee
1}\big]\},$$ and we conclude by Theorem \ref{GK}.

\hfill$\Box$

\begin{remark} \rm As clearly stated at the beginning of this section, Theorem \ref{princ} holds in the case $\alpha >1$. We believe that this is due to the particular arguments used for the proof, and that it is possible to extend the ASLLT also to the case $\alpha<1$. The critical case $\alpha=1$ remains unexplored till now. Another not yet investigated situation is
for $\alpha=2$ with $x \mapsto E[X^2 1_{\{|X|\leqslant x\}}]$ slowly varying and $E[X^2]=\infty$ with $x \mapsto x^2P(|X| > x)$ not
slowly varying. Hopefully, we shall treat these cases in another paper.
\end{remark}

\end{document}